\documentclass{gtart_h}

\def\ifplaintex{\expandafter\ifx\csname documentclass\endcsname\relax}

\def\gtp{{\mathsurround=0pt\it $\cal G\mskip-2mu$eometry \&\ 
$\cal T\!\!$opology $\cal P\!$ublications}}  

\def\recd{{\small Received:\qua\receiveddate\ifx\reviseddate\relax
\else\qquad Revised:\qua\reviseddate\fi\par}} 

\def\lognumber#1{\def\thelognumber{#1}}
\def\volumenumber#1{\def\thevolumenumber{#1}}
\def\volumeyear#1{\def\thevolumeyear{#1}}
\def\papernumber#1{\def\thepapernumber{#1}}
\def\pagenumbers#1#2{\def\startpage{#1}\def\finishpage{#2}}
\def\published#1{\def\publishdate{#1}}

\def\received#1{\def\receiveddate{#1}}

\def\accepted#1{\def\accepteddate{#1}}

\def\asciiemail#1{\def\theasciiemail{#1}}
\def\asciiurl#1{\def\theasciiurl{#1}}

\long\def\asciiabstract#1{\long\def\theasciiabstract{#1}}

\let\\\par\let\thelognumber\relax\let\thevolumenumber\relax
\let\thepapernumber\relax\let\thevolumeyear\relax\let\startpage\relax
\let\finishpage\relax\let\publishdate\relax\let\receiveddate\relax
\let\reviseddate\relax\let\accepteddate\relax\let\theasciititle\relax
\let\theasciiauthors\relax
\let\theasciiabstract\relax

\let\theasciiemail\relax
\let\theasciiurl\relax

\ifplaintex
\font\logobig=cmssbx10 scaled 3836
\font\logomed=cmssbx10 scaled 2557
\else
\font\logobig=cmssbx10 scaled 4200
\font\logomed=cmssbx10 scaled 2800
\fi

\long\def\makeagttitle{   
\count0=\startpage
\agt\hfill      
\hbox to 45truept{\vbox to 0pt{\vglue -13truept{\logomed A\kern -.37em{\logobig 
T}\kern -.38em G}\vss}\hss}
\break
{\small Volume \thevolumenumber\ (\thevolumeyear)
\startpage--\finishpage\nl
Published: \publishdate}

\vglue .25truein

{\parskip=0pt\leftskip 0pt plus
1fil\def\\{\par\smallskip}{\Large\bf\thetitle}\par\medskip} \vglue
0.05truein

{\parskip=0pt\leftskip 0pt plus 1fil\def\\{\par}{\sc\theauthors}
\par\medskip}%
 
\vglue 0.03truein 

{\small\leftskip 25truept\rightskip 25truept{\bf Abstract}\stdspace\theabstract

{\bf AMS Classification}\stdspace\theprimaryclass
\ifx\thesecondaryclass\relax\else; \thesecondaryclass\fi\par
{\bf Keywords}\stdspace \thekeywords\par}\vglue 7truept

}   

\ifplaintex
\font\phead=cmsl9 scaled 950
\font\pnum=cmbx10 scaled 913
\font\pfoot=cmsl9 scaled 950
\headline{\vbox to 0pt{\vskip -4.5mm\line{\small\phead\ifnum
\count0=\startpage ISSN 1472-2739 (on-line) 1472-2747 (printed)
\hfill {\pnum\folio}\else\ifodd\count0\def\\{ }%
\ifx\theshorttitle\relax\thetitle\else\theshorttitle\fi\hfill{\pnum\folio}
\else\def\\{ and }{\pnum\folio}\hfill\ifx\theshortauthors\relax\theauthors
\else\theshortauthors\fi\fi\fi}\vss}}
\footline{\vbox to 0pt{\vglue 0mm\line{\small\pfoot\ifnum\count0=\startpage
\copyright\ \gtp\hfill\else
\agt, Volume \thevolumenumber\ (\thevolumeyear)\hfill\fi}\vss}}
\else
\font=cmsl9 scaled 1050
\font\lnum=cmbx10 
\font=cmsl9 scaled 1050
\makeatletter
\def\@oddhead{{\small\lhead\ifnum\count0=\startpage ISSN 1472-2739 
(on-line) 1472-2747 (printed)\hfill {\lnum\number\count0}\else\ifodd\count0
\def\\{ }\ifx\theshorttitle\relax \thetitle \else\theshorttitle\fi\hfill
{\lnum\number\count0}\else\def\\{ and }{\lnum\number\count0}
\hfill\ifx\theshortauthors\relax 
\theauthors\else\theshortauthors\fi\fi\fi}}\def\@evenhead{\@oddhead}
\def\@oddfoot{\small\lfoot\ifnum\count0=\startpage\copyright\ \gtp\hfill\else
\agt, Volume \thevolumenumber\ (\thevolumeyear)\hfill\fi}
\def\@evenfoot{\@oddfoot}
\makeatother
\fi
\let\maketitlepage\makeagttitle
\let\makeshorttitle\maketitlepage
\let\maketitle\maketitlepage

\newwrite\gtoutfile
\long\gdef\makeheadfile{  
{\def\\{, }\def\s{ }
\immediate\openout\gtoutfile head.xxx
\immediate\write\gtoutfile{Proxy-for: \ifx\theasciiauthors\relax
\theauthors\else\theasciiauthors\fi\s<\ifx\theasciiemail\relax\theemail\else\theasciiemail\fi>}
\immediate\write\gtoutfile{\noexpand\\}
\immediate\write\gtoutfile{Authors: \ifx\theasciiauthors\relax
\theauthors\else\theasciiauthors\fi}
{\def\\{ }\immediate\write\gtoutfile{Title: \ifx\theasciititle\relax
\thetitle\else\theasciititle\fi}}
\immediate\write\gtoutfile{Subj-class: GT or SG, GR etc}
\immediate\write\gtoutfile{MSC-class: \theprimaryclass\ifx\thesecondaryclass\relax\else, \thesecondaryclass\fi}
\immediate\write\gtoutfile{Journal-ref: Algebr. Geom. Topol. \thevolumenumber\s
(\thevolumeyear) \startpage-\finishpage}
\immediate\write\gtoutfile{Comments: Published by Algebraic and
Geometric Topology at}
\immediate\write\gtoutfile{\s\s\s  http://www.maths.warwick.ac.uk/agt/AGTVol\thevolumenumber/agt-\thevolumenumber-\thepapernumber.abs.html}
\immediate\write\gtoutfile{\noexpand\\}
\immediate\write\gtoutfile{}
\ifx\theasciiabstract\relax
\immediate\write\gtoutfile{\theabstract}\else
\immediate\write\gtoutfile{\theasciiabstract}\fi
\immediate\write\gtoutfile{}
\immediate\write\gtoutfile{\noexpand\\}
\immediate\write\gtoutfile{}
\immediate\closeout\gtoutfile}}  

\def\maketitlepage{\makeagttitle\makeheadfile}
\let\makeshorttitle\maketitlepage
\let\maketitle\maketitlepage

\lognumber{48}
\volumenumber{4}
\volumeyear{2004}
\papernumber{48}
\pagenumbers{1111}{1123}
\received{22 September 2004} 
\accepted{15 November 2004}
\published{25 November 2004}

\usepackage{amsmath,amssymb}

\newtheorem{theo}{Theorem}[section]

\newtheorem{lem}[theo]{Lemma}
\newtheorem{defi}[theo]{Definition}

\newtheorem{corol}[theo]{Corollary}

\def \Z{{\bf Z}}

\def \s{$\Sigma$}

\def\nt{$n$-trivializer}

\begin{document}

\title[Alexander polynomial and hyperbolic volume]{Alexander
polynomial, finite type invariants\\and volume of hyperbolic knots}

\author{Efstratia Kalfagianni}

\address {Department of Mathematics, Michigan State
University, E. Lansing, MI 48824, USA\\or\\School of Mathematics, 
Institute for Advanced Study, Princeton, NJ 08540, USA}

\gtemail{\mailto{kalfagia@math.msu.edu}{\rm\qua 
or\qua}\mailto{kalfagia@math.ias.edu}}
\asciiemail{kalfagia@math.msu.edu, kalfagia@math.ias.edu}

\urladdr{http://www.math.msu.edu/~kalfagia}
\asciiurl{http://www.math.msu.edu/ kalfagia}

\begin{abstract} We
show that given $n>0$, there exists a hyperbolic knot $K$ with trivial Alexander
polynomial, trivial finite type invariants of order $\leq n$, and such that the volume of the complement
of $K$ is larger than $n$. This contrasts with the
known statement that the volume of the complement of a hyperbolic alternating knot is bounded above by a linear function of the coefficients of the Alexander polynomial of the knot.
As a corollary to our main result we obtain
that, for every $m>0$, there exists a sequence of hyperbolic  knots
with trivial
finite type invariants of order $\leq m$ but arbitrarily large volume.
We discuss how our results fit within the framework of relations
between the finite type invariants and the volume of hyperbolic
knots, predicted by Kashaev's hyperbolic volume conjecture.
\end{abstract}

\asciiabstract{%
We show that given n>0, there exists a hyperbolic knot K with trivial
Alexander polynomial, trivial finite type invariants of order <=n, and
such that the volume of the complement of K is larger than n. This
contrasts with the known statement that the volume of the complement
of a hyperbolic alternating knot is bounded above by a linear function
of the coefficients of the Alexander polynomial of the knot.  As a
corollary to our main result we obtain that, for every m>0, there
exists a sequence of hyperbolic knots with trivial finite type
invariants of order <=m but arbitrarily large volume.  We discuss how
our results fit within the framework of relations between the finite
type invariants and the volume of hyperbolic knots, predicted by
Kashaev's hyperbolic volume conjecture.}



\maketitle

\section{Introduction}
Let $c(K)$ denote the crossing number
and let
$\Delta_K(t):= \sum_{i=0}^k c_it^i$
denote the Alexander  polynomial
of a knot $K$. If $K$ is hyperbolic, let ${\rm vol}(S^3\setminus K)$ denote the volume of its complement.
The determinant of $K$ 
is the quantity ${\rm det}(K):= |\Delta_K(-1)|$. Thus, in general,
we have
$${\rm det}(K)\leq ||\Delta_K(t)||:=\sum_{i=0}^k|c_i|.\eqno(1)$$
It is well know that
the degree
of the Alexander polynomial of an alternating knot equals twice the genus of the knot.
On the other hand one can easily construct knots with arbitrarily large genus
and trivial Alexander polynomial by taking connected sums of untwisted
Whitehead doubles. In this note, we construct hyperbolic knots of arbitrarily large
genus and volume that have trivial Alexander polynomial. To put our results into the
appropriate 
context, let us first recall what is known about the relation between volume and Alexander polynomial of
alternating knots.
In general, it is known (see \cite{car}
and references therein) that there
exists a universal constant $C>0$ such that
if $K$
is  a hyperbolic knot with crossing number $c(K)$, then
$${\rm vol}(S^3\setminus K)\leq C c(K). \eqno(2)$$
If in addition $K$ is alternating then one also has ( \cite{bz})
$$c(K)\leq {\rm det}(K)\ \ {\rm and}\ \ {\rm det}(K)=||\Delta_K(t)||.
\eqno(3)$$
By a result of Menasco \cite{me},
prime, alternating knots are either torus  knots or hyperbolic.
Combining this with (1)-(3),  we derive that
there is a universal constant $C>0$ such that we have
$$ {\rm vol}(S^3\setminus K)\leq C {\rm det}(K)\ \ {\rm and} \ \ {\rm vol}(S^3\setminus K)\leq C||\Delta_K(t)||,
 \eqno(4)
 $$
for every prime, alternating non-torus knot $K$. (In fact, in \cite{sto}, A. Stoimenow
showed that the quantity ${\rm det}(K)$ in (4) can be replaced by $\log ({\rm det}(K))$.)
Thus, the Alexander polynomial of prime alternating knots, dominates two of the most important geometric knot invariants; namely the volume and the genus. In contrast to that,
the main result in this paper implies
the following.
 \begin{corol} \label{corol:Alexander} Given $n\in {\bf N}$, there exists a hyperbolic knot $K$
with trivial Alexander polynomial for which we have
${\rm vol }(S^3\setminus K)>n$ and ${{\rm genus}(K) > {\textstyle {n}\over 6}} $.
\end{corol}
The development of quantum topology
led to many new knot invariants that generalized  the Alexander 
polynomial. These are the Jones type
polynomial invariants (quantum invariants)  and their perturbative counterparts known as
finite type invariants (or Vassiliev invariants). An open conjecture relating these invariants to the geometry
of the knot complement is the hyperbolic volume conjecture
of R. Kashaev \cite{ka}.
Kashaev's conjecture, as rephrased after the
work of H. Murakami and J. Murakami \cite{mu},
asserts that the volume of the complement of a hyperbolic knot 
$K$ is equal to a certain limit of values of the  colored Jones polynomials of $K$.
A convenient way to organize these  polynomials is via the so-called 
colored Jones function. This is a 2-variable formal power series
$$ J_K(h, d)= \sum_{i\geq 0, j\leq i} a_{ij} h^i d^j  \ \in {\bf Q}[h, d], $$
\noindent  such that $a_{ij}$ is a Vassiliev invariant of order $i$ for $K$.
By the Melvin-Morton-Rozansky conjecture, the first proof of which was given in
\cite{bng},  the diagonal
subseries $D_K(h):={ \sum_{i} a_{ii} h^i}$ is essentially equivalent to the Alexander
polynomial of $K$.
The volume conjecture  implies that the colored Jones function
of a hyperbolic knot $K$ determines the volume of the complement of $K$.
In particular, it implies that if there exist two hyperbolic knots that
have {\em all} of their finite type invariants the same then their
complements should have equal volumes.
In this setting, Corollary \ref{corol:Alexander}
says that the part $D_K(h)$ of $ J_K(h, d)$ is far from controlling the volume
of hyperbolic knots. The following theorem, which is the main result of the paper,
provides a stronger assertion. It implies that there can be no universal value
$n\in {\bf N}$ such that, for every hyperbolic knot $K$,
the finite type invariants of orders $\leq n$
together with $D_K(h)$
determine the volume of the complement of $K$.
\begin{theo} \label{theo:main} For every $n\in {\bf N}$, there exists a hyperbolic knot $K$ such that:

\noindent {\rm(a)}\qua We have $\Delta_K(t)=1$.

\noindent {\rm(b)}\qua All the finite type invariants of orders $<2n-1$ vanish for $K$.

\noindent {\rm(c)}\qua We have ${\rm vol }(S^3\setminus K)>n$.

\noindent {\rm(d)}\qua We have ${\rm genus }( K)>{\textstyle {n}\over 6}$.
\end{theo}
Theorem \ref{theo:main} has the following corollary
which implies that for every $m\in {\bf N}$, there are hyperbolic knots with
arbitrarily large volume but
trivial ``$m$- truncated" colored Jones function. The proof of the corollary is spelled out at the end of Section 3.
\begin{corol} \label{corol:infinite} For every $n, m\in {\bf N}$,
there exists a hyperbolic knot $K$ with trivial finite type invariants of orders $\leq m$
and such that
${\rm vol }(S^3\setminus K)>n$.
\end{corol}
The proof of Theorem \ref{theo:main} combines  results
of \cite{ak} and \cite{kl1} with
techniques from hyperbolic geometry \cite{thurston}.
Given $n\in {\bf N}$, \cite{ak} provides a method for constructing knots with
trivial Alexander polynomial and trivial finite type invariants of orders $<2n-1$.
To obtain hyperbolicity and appropriate volume growth we combine a result from
\cite{kl1} with the work of Thurston and a result of Adams \cite{adams}.
The results of \cite{ak}, \cite{kl1} needed in this paper are summarized in Section 2
where we also prove some auxiliary lemmas.
Section 3 is devoted to the proof of Theorem \ref{theo:main}.

\medskip
{\bf Acknowledgment}\qua The author thanks Oliver Dasbach and
Xiao-Song Lin for useful discussions and correspondence.  She also
thanks Alex Stoimenow for his interest in this work and Dan Silver for
helpful comments and for informing her about his work \cite{sstw}. In
addition, she acknowledges the support of the NSF through grants
DMS-0104000 and DMS-0306995 and of the Institute for Advanced Study
through a research grant.

\section {Preliminaries} In this section we recall results and terminology
that will be used in the subsequent sections and we prove some auxiliary lemmas needed for the proof of the main
results.
A crossing disc for a knot $K$ is a disc $D$ that intersects $K$ only in its interior
exactly twice with zero algebraic intersection number; the boundary $\partial D$ is called a crossing circle.
Let $q\in \Z$.
A knot $K'$ obtained from $K$ by twisting  $q$-times along $D$
is said to be obtained by a generalized crossing change of order $|q|$. Clearly,
if $q=0$ we have $K=K'$ and if $|q|=1$ then $K,K'$ differ by an ordinary crossing change.
Note that
$K'$ can also be viewed as the result from $K$ under $\textstyle{1\over q}$-surgery of $S^3$
along $\partial D$. 
\begin{defi} \label{defi:adj}{\rm \cite{kl1}}\qua
We
will say that $K$ is $n$-adjacent to the unknot, for some
$n\in {\bf N}$,
if $K$ admits an embedding containing $n$ generalized crossings
such that changing any $0<m\leq n$ of them yields an
embedding of  the unknot. A collection of crossing circles
corresponding to these crossings is called an \nt.
\end{defi}

Let $V$ be a solid torus in $S^3$ and suppose that a
knot $K$ is embedded in $V$ so that it is geometrically essential. We will use the term
``$K$ is $n$-adjacent to the unknot in $V$" to mean the following:
There exists an embedding of $K$ in $V$
that contains $n$ generalized crossings
such that changing any $0<m\leq n$ of them transforms $K$
into a knot that bounds an embedded disc in {\rm int}V.

Recall that if $K$ is
a non-trivial satellite with
{\em companion knot} ${\hat C}$  and {\em model knot} $\hat K$ then:

i)\qua ${\hat C}$ is non-trivial;

ii)\qua $\hat K$ is geometrically essential in a standardly embedded
solid torus $V_1\subset S^3$ but not isotopic to the core of $V_1$;  and

iii)\qua there is a homeomorphism
$h: V_1\longrightarrow  V:=h(V_1),$
such that $h(\hat K)=K$ and $\hat C$ is the core of $V$.

The following theorem will be used to ensure that for every $n>0$ there is a hyperbolic knot that is $n$-adjacent to the unknot.
\begin{theo} \label{theo:satellite} {\rm (\cite{kl1}, Theorem 1.3)} Let $K$ be a non-trivial satellite knot
and let $V$ be any companion solid torus of $K$.
If
$K$ is
$n$-adjacent to the unknot, for some $n>0$, then
it is $n$-adjacent to the unknot in $V$. Furthermore,
any model knot of $K$ is 
$n$-adjacent to the unknot in the standard solid torus $V_1$.
\end{theo}

The next result compiles the statements of Theorems 4.5, 5.2 and 5.6 of \cite{ak}.

\begin{theo} \label{theo:invariants} {\rm \cite{ak}}\qua
Suppose that $K$ is a knot that is  $n$-adjacent to the unknot
for some $n>2$. Then, the following are true:

\noindent {\rm(a)}\qua  All finite type invariants
of
orders less than $2n-1$ vanish for $K$.

\noindent {\rm(b)}\qua We have $\Delta_K(t)=1$.
\end{theo}

Given a knot $K$, a collection of crossing circles $K_1, \ldots, K_n$
and an $n$-tuple of integers ${\bf r}:=(r_1 \ldots, r_n)$,
let $K( {\bf r})$ denote the knot obtained from $K$ by performing a generalized crossing of order $|r_i|$ along $K_i$,
for all $1\leq i\leq n$.
\begin{lem} \label{lem:twisting} Suppose that $K$ is $n$-adjacent to the unknot, for some $n>1$,
and that the crossing circles $K_1, \ldots, K_n$ constitute  an \nt \ for $K$. Then,
for every $\bf r$ as above, 
$K({\bf r})$ is also $n$-adjacent to the unknot.
\end{lem}
\proof 
It follows immediately from Theorem 4.4 and Remark 5.5 of \cite{ak}. \endproof

For a knot $K$
that is $n$-adjacent to the unknot and an \nt \
$L_n:=K_1\cup \ldots \cup K_n$ let
$\eta{(L_n\cup K)}$ denote a
tubular neighborhood of $L_n\cup K$. We close this section
with two technical lemmas that we need in the next section.

\begin{lem} \label{lem:irreducible}
Suppose that $K$ is a non-trivial knot
that is $n$-adjacent to the unknot, for some $n>0$. Then for every \nt\ $L_n$,
the 3-manifold  $M(K, L_n):= S^3\setminus \eta{(L_n\cup K)}$ is irreducible and
$\partial$-incompressible.
\end{lem}

\proof Let $D_1,\ldots, D_n$ be crossing discs bounded by
$K_1,\ldots, K_n$, respectively.
Assume, on the contrary, that $M(K, L_n)$
contains an essential 2-sphere $\Sigma$.
Assume that $\Sigma$ has been isotoped so that the intersection
$I:=\Sigma \cap ({\cup_{i=1}^nD_i})$ is minimal.
We must have $I\neq \emptyset$
since otherwise $\Sigma$ would bound
a 3-ball in $M(K, L_n)$. Let $c\in (\Sigma\cap D_i)$ denote a component of $I$ that is innermost
on $\Sigma$; it bounds a disc $E\subset \Sigma$ whose interior is
disjoint from
${\cup_{i=1}^n D_i}$.
Since $\Sigma$ is separating in $M(K, L_n)$, $E$ can't contain just one point of $K\cap D_i$
and, by the minimality assumption,
$E$ can't be disjoint from $K$. Hence $E$ contains both points of $K\cap D_i$ and so $c=\partial E$ is parallel to $\partial D_i$ in $D_i\setminus K$.
It follows that $K_i$ bounds an embedded disc
in the complement of $K$. But then
no twist along $K_i$ can unknot $K$, contrary to
our assumption that $K_i$ is part of an \nt \ of $K$.
This contradiction
finishes the irreducibility claim.
To finish the proof of the lemma, suppose that $M(K, L_n)$
is $\partial$-compressible. Then, a component $T$ of $\partial M(K, L_n)$
admits a compressing disc $E'$. Since $K$ is non-trivial, we must have
$T=\partial \eta(K_i)$, for some $1\leq i\leq n$. But then, $\partial E$ must be a longitude of $\eta(K_i)$. Thus $K_i$ bounds a disc in the complement of $K$; a contradiction. \endproof

\begin{lem} \label{lem:nonisotopic} Suppose that $K$ is a knot that is embedded in
a standard solid torus $V_1\subset S^3$ so that it is geometrically essential in $V_1$.  Suppose, moreover,
that $K$ is $n$-adjacent to the unknot in $V_1$, for some $n>1$.
Then, $K$ is not the unknot.
\end{lem}
\proof Let $L_n$ be an  \nt \  that exhibits $K$ as $n$-adjacent to the unknot in $V_1$.
Let $\cal D$ denote the disjoint union of a collection of crossing discs corresponding to the
components of $L_n$. Suppose, on the contrary, that $K$ is the unknot.
By Lemma 4.1 and the discussion in Remark 5.5 of \cite{ak}, $K$ bounds an
embedded disc $\Delta$
in the complement of $L_n$.
Furthermore,  $\Delta \cap {\cal D}$ is a
collection $\cal A$ of $n$ disjoint properly embedded arcs in $\Delta$;
one for each component of $\cal D$. Now making a
crossing change on $K$ supported on a component  $K_i\subset L_n$
is the same as twisting $\Delta$ along the arc in $\cal A$
corresponding to $K_i$. Since $n>1$, $\cal A$ has at least two components, say
$\alpha_1, \alpha_2$, which cut $\Delta$ into three subdiscs. Since $K$ is geometrically essential
in $V$, but can be unknotted by twisting along either of  $\alpha_1$, $\alpha_2$,
exactly
two of these subdiscs must intersect $\partial V_1$ in a longitude.
Recall that $\alpha_1, \alpha_2$
correspond to crossings that exhibit $K$ as 2-adjacent to the unknot in $V_1$.
Thus $K$ can be also be uknotted by
performing
the aforementioned twists simultaneously on $\alpha_1\cup \alpha_2$.
A straightforward checking will convince the reader that this is impossible. \endproof

\section{Knot adjacency and hyperbolic volume}
\subsection{Trivializers and hyperbolic Dehn filling }
In this section we will prove the main results of this note.
As we will show
Theorem \ref{theo:main} follows from the following theorem
and known hyperbolic volume techniques.
\begin{theo} \label{theo:hyperbolic} For
every $n>2$, there exists a knot ${\hat K}$ that is $n$-adjacent to the unknot and it
admits an \nt \ $L_n$ such that the interior of
$M_n:=M({\hat K}, L_n)$ admits a complete
hyperbolic structure with finite volume. Furthermore,
we have  ${\rm vol}(M_n)>n\,v_3 $, where
$v_3\, (\approx 1.01494)$ is the volume of a regular hyperbolic ideal tetrahedron.
\end{theo}

The proof of Theorem \ref{theo:hyperbolic} uses Thurston's uniformization theorem for Haken 3-manifolds
\cite{thurston}
in combination with a result of Adams \cite{adams}.
Having Theorem \ref{theo:hyperbolic} at hand, Theorem \ref{theo:main} follows from
Thurston's hyperbolic Dehn filling theorem. The
statements of these results of Thurston can also be found in \cite{bo}.
By abusing the terminology, we will say $M_n$ is hyperbolic
to mean that the interior of $M_n$ admits a complete
hyperbolic structure with finite volume.

Next let us show how Theorem \ref{theo:main} follows from Theorem \ref{theo:hyperbolic}.
For the convenience of the reader we repeat the statement of Theorem \ref{theo:main}.

\begin{theo} For every $n\in {\bf N}$, there exists a hyperbolic knot $K$ such that:

\noindent {\rm(a)}\qua We have $\Delta_K(t)=1$.

\noindent {\rm(b)}\qua All the finite type invariants of orders $<2n-1$ vanish for $K$.

\noindent {\rm(c)}\qua We have ${\rm vol }(S^3\setminus K)>n$.

\noindent {\rm(d)}\qua We have ${\rm genus }( K)>{\textstyle {n}\over 6}$.
\end{theo}

\proof  First assume
that $n\leq 2$. Then, we have  $2n-1\leq 3$ and thus to satisfy part(c)  we need a knot that has trivial finite type invariants of orders $\leq 2$.
It is known that the only finite type invariant of order $\leq 2$
is essentially the  second derivative of the Alexander polynomial
at $t=1$. Thus, a knot with trivial Alexander polynomial
has trivial finite type invariants of order $\leq 2$. There are hyperbolic
knots with trivial Alexander polynomial;
the first such knots occur at thirteen crossings and all  have volume $>2$.
Thus for $n\leq 2$ there are hyperbolic knots satisfying parts (a)-(c).
But note that (d) is trivially satisfied.

Suppose now that $n>2$.
By Theorem \ref{theo:hyperbolic}, there exists a non-trivial knot
${\hat K}$, that is $n$-adjacent to the unknot,
and an \nt, \ say $L_n$, so that
$M_n$ is hyperbolic and ${\rm vol}(M_n)> n\ v_3 $.
Let $K_1, \ldots, K_n$ denote the components of $L_n$. Given an
$n$-tuple of  integers ${\bf r}:= (r_1, \ldots, r_n)$,
let $M_n({\bf r})$ denote the 3-manifold obtained
from $M_n$ as follows: For $1\leq i\leq n$, perform
Dehn filling with slope ${\textstyle {1}\over {r_i}}$-surgery along $\partial {\eta(K_i)}$.
Let ${\hat K}({\bf r})$ denote the image of $\hat K$ in $M_n({\bf r})$.
Clearly, $M_n({\bf r}):=S^3\setminus \eta({\hat K}({\bf r}))$ and
${\hat K}({\bf r})$ is obtained from $\hat K$ by generalized crossing
changes.
By Lemma \ref{lem:twisting}, ${\hat K}({\bf r})$ is $n$-adjacent to
the unknot, for every $n$-tuple ${\bf r}$. On the other hand,
by
Thurston's hyperbolic Dehn filling theorem, if $r_i>> 0$
then $M_n({\bf r})$ admits a complete hyperbolic structure
of finite volume; thus ${\hat K}({\bf r})$ is a hyperbolic knot.
By the proof of Thurston's theorem, the hyperbolic metric on
$M_n({\bf r})$ can be chosen so that it is arbitrarily close
to the metric of $M_n$, provided that the numbers
$r_i$ are all sufficiently large.  Thus by choosing
the $r_i$'s large
we may ensure that the volume of $M_n({\bf r})$
is arbitrarily close to that of $M_n$. 
Since ${\rm vol}(M_n)>n\ v_3 $,
we can choose the $r_i$'s so that we have
${\rm vol}(M_n({\bf r}))> n$. But
since ${\rm vol}(M_n({\bf r}))={\rm vol}(S^3\setminus {\hat K}({\bf r}))$, 
setting $K:={\hat K}({\bf r})$
we are done. This finishes the proof of parts (a)-(c) of the theorem.
For part (d), recall that by the main result of \cite{hlu}
(or by Theorem 1.3 of \cite{kl2}), we have
$6 {\rm genus}(K)-3 \geq n$. Thus,
${\rm genus}(K)\geq {\textstyle {n+3}\over {6}}> {\textstyle {n}\over {6}}$
as desired.
\qed

\subsection{ Ensuring hyperbolicity}
The rest of the section is devoted to the proof of Theorem \ref{theo:hyperbolic}.
Our purpose is to show that given  $n>2$,  we can find a non-trivial knot
 $\hat K$, that is exhibited $n$-adjacent to the unknot by an \nt \ $L_n$ such that
 the 3-manifold
$M_n:=S^3\setminus \eta ({\hat K}\cup L_n)$ has the following properties:

(i)\qua $M_n$ is irreducible.

(ii)\qua $M_n$ is $\partial$-incompressible.

(iii)\qua $M_n$ is atoroidal; every incompressible embedded torus is parallel to a component of $\partial M_n$.

(iv)\qua $M_n$ is anannular; every incompressible  properly embedded annulus
can be isotoped on $\partial M_n$.

Having properties (i)-(iv) at hand, we can apply Thurston's
uniformization theorem for Haken 3-manifolds  to conclude
that $M_n$
is hyperbolic. Now $\partial M_n$ is a collection of tori.
Associated with each component of $\partial M_n$,
there is a cusp in ${\rm int}M_n$ homeomorphic to $T^2\times [1, \infty)$.
A result of Adams \cite{adams} states
that a complete hyperbolic 3-manifold with $N$ cusps must have
volume at least $N. \ v_3$. Since  $\partial M_n$
has $n+1$ components, ${\rm int}M_n$
has $n+1$ cusps. Thus
${\rm vol}(M_n)\geq (n+1)\ v_3 > n\ v_3$ and the desired conclusion follows.

By \cite{ak}, for every  $n\in {\bf N}$, there exist many non-trivial knots that are
$n$-adjacent to the unknot. Fix
$n>2$.  Let $\bar K$ be any non-trivial knot that is $n$-adjacent to the unknot
and let $L_n$ be any \nt \ corresponding to $\bar K$.
By lemma \ref{lem:irreducible},
properties (i)-(ii) hold for $M_n:= S^3\setminus \eta ({\bar K}\cup L_n)$.
Next we set off to show that we can  choose $\bar K$ and $L_n$ so that
$M_n$ satisfies properties
(iii)-(iv) as well. We need the following lemma that will be used
to control essential tori in the complement of \nt's.
\begin{lem} \label{lem:zero} Let $L_n$ be an \nt \ that shows a non-trivial knot $\bar K$ to be $n$-adjacent to the unknot.
Suppose that $M_n:=S^3\setminus \eta(L_n\cup {\bar K})$ contains an essential torus $T$
and let $V$ denote the solid torus in $S^3$ bounded by $T$. 
Then, $L_n$ can be isotoped in ${\rm int} V$
by an isotopy
in the complement of $\bar K$.
Furthermore, $T$ remains essential in the complement of $\bar K$.
\end{lem}
\proof First we show that $L_n$ can be isotoped
in ${\rm int} V$, by an isotopy
in the complement of $\bar K$.
Suppose that a component $K_1\subset L_n$ lies outside $V$
and let
$D_1$ be a crossing disc bounded by $K_1$. We can isotope $D_1$ so that
every component of $D_1\cap T$ is essential in $D_1\setminus {\bar K}$. Thus, each component of 
$D_1\cap T$ is either parallel to $\partial D_1$ on $D_1$, or it bounds a
subdisc of $D_1$ that is intersected exactly once by $\bar K$. In the first case, $K_1$ can be
isotoped to lie in ${\rm int}V$. In the later case,  each component of
$D_1\cap T$  bounds a subdisc $E\subset D_1$ that is intersected exactly once by $\bar K$.
Since $T$ is essential in the complement of $\bar K$, $E$ must be a meridian disc of $V$.
Thus $\bar K$ has wrapping number one in the follow-swallow torus. It follows that 
$\bar K$ is composite and that
the crossing change corresponding to $L_1$ occurs within a factor of $\bar K$.
But then such a crossing change cannot unknot $\bar K$ contradicting the fact
that $L_1$ is part of an \nt. Thus
$L_n$ can be isotoped to lie in ${\rm int} V$ as desired. Now it is easy to see
that $L_n$
exhibits ${\bar  K}$ as $n$-adjacent to the unknot in ${\rm int} V$.

Next we show that $T$ remains essential in the complement
of $\bar K$.

Since the linking number of $\bar K$ and each component of $L_n$
is zero, $\bar K$ bounds a Seifert surface in the complement of $L_n$.
Let $S$ be a Seifert surface of $\bar K$ that is of minimum genus
(and thus of maximum Euler characteristic) among all
such Seifert
surfaces of $\bar K$. Since $T$, $S$ are incompressible in $M_n$
we can isotope them so that
$S\cap T$ is a collection of parallel essential curves on $T$.
Since twisting along any component of $L_n$ unknots $\bar K$,
the winding number ${\bar K}$ in $V$
is zero.
We conclude that $S\cap T$ is homologically trivial in $T$.
Thus we may
replace the components of $S\setminus T$  that lie outside $V$ by annuli
on $T$, and then isotope these annuli away from $T$ in $V$, to obtained a Seifert surface $S'$
 for
$\bar K$. Now $S'$ lies in $M_n\cap V$. Since no component
of $X:=S\cap {\overline {S^3\setminus V}}$ can be a
disc, the Euler characteristic $\chi(X)$ is not positive.
Thus $\chi(S')\geq \chi(S)$. Since $S$ was chosen to be of maximum
Euler characteristic for $K$ in the complement of $L_n$,
we have $\chi(S')\leq \chi(S)$.
It follows
that $\chi(S')= \chi(S)$ and thus
${\rm genus}(S')={\rm genus}(S)$. Hence $S'$ is also a minimum genus 
Seifert surface for $\bar K$ in $M_n$.
Let $\cal D$  be a collection
of crossing discs; one for each component of $L_n$.
For $D\in {\cal D}$,
the components of $S'\cap D$ are closed curves that are parallel to $\partial D$
on $D$ and an arc properly embedded on $S'$. After an isotopy of $S'$ we can arrange so that
$S'\cap D$ is only an arc $\alpha$ properly embedded in $S'$.
Thus, $S' \cap {\cal D}$ is a collection, say ${\cal A}$
of arcs properly embedded on $S'$. 
By Theorem 4.1 of \cite{hlu}, $S'$ remains of minimum genus
in the complement of  $\bar K$. That is, ${\rm genus}({\bar K})={\rm genus}(S')$.
Thus $S'$ is incompressible in the complement of
$\bar K$.
By assumption, twisting along all the components of $L_n$ turns $\bar K$ into the unknot in $V$.
This unknot can be isotoped to lie in a 3-ball $B\subset {\rm int}V$.
Assume that $\bar K$ is inessential in $V$; then
$\bar K$ can also be isotoped to lie in  $B$.
Since $S'$ is incompressible, $S'$ can  also be isotoped to lie in $B$. But
then, $\cal A$ and $\bar K$ will lie in $B$.
Now since $L_n$ can be isotoped
so that it lies in a small neighborhood of ${\cal A}$
in a collar of $S'$,
it follows that ${\bar K}\cup L_n$ can be isotoped in $B$.
But this contradicts the assumption that $T$ is essential in $M_n$. \endproof

Next we turn our attention to essential annuli in $M_n$.
It is known that an atoroidal link complement that contains essential annuli admits
a Seifert fibration over a sphere with at most three punctures.
In particular such a link can have at most 3-components. Hence, the restriction
$n>2$ in our case implies that if $M_n$ is atoroidal then it is anannular.
For the convenience of the reader we give a direct argument that $M_n$ is anannular
in the next lemma.
\begin{lem} \label{lem:annuli} Let $L_n$ be an \nt \  of a non-trivial knot  $\bar K$.
If $n>2$ and $M_n$ is atoroidal, then $M_n$ is anannular.
\end{lem}
\proof
Since $\bar K$ is non-trivial, by Lemma \ref{lem:irreducible}, $M_n$ is irreducible.
We will show that every incompressible, properly embedded annulus $(A, \partial A)\hookrightarrow (M_n, \partial M_n)$
can be isotoped to lie on $\partial M_n$.

First, suppose that $A$ runs between two components of $\partial M_n$, say $T_1, T_2$.
Let $N$ denote a neighborhood of $A\cup T_1\cup T_2$. Now $\partial N$ has three tori components.
Two of them are parallel to $T_1, T_2$, respectively. Let $T$ denote the third one. Since $n>2$,
$\partial M_n$ has at least four components. Thus  $T$ separates pairs of components
of $\partial M_n$ and it can't be boundary parallel. Hence, since $M_n$ is atoroidal,
$T$ must be compressible. A compressing disc cuts $T$ into a 2-sphere that
separates pairs of components
of $\partial M_n$ and thus it is essential; a contradiction.
Therefore, $\partial A$ must lie on one component, say $T_1$, of $\partial M_n$.
Let $N$ denote a neighborhood on $A\cup T_1$. Now $\partial N$ contains
a torus $T$ that separates $T_1$ from the rest of the components of $\partial M_n$.
Suppose, for a moment, that $T$ compresses. Then a compressing disc cuts $T$ into a 2-sphere that separates
$T_1$ from the rest of the components of $\partial M_n$. But
such a sphere would be essential contradicting the irreducibility of $M_n$.
Therefore we conclude that $T$ must be parallel to a component of $\partial M_n$.
But since $T$ separates $T_1$ from at least three components of $\partial M_n$,
it can only be 
parallel to $T_1$. Hence $A$ is contained in a collar $T_1\times I$
and it can be isotoped on $T_1$. This proves that $A$ is inessential which finishes the proof of the lemma. \endproof

\subsection{Completing the proofs  }
We are now ready to give the proofs of the main results.
We begin with the proof of Theorem \ref{theo:hyperbolic}.

\begin{theo}  For
every $n>2$, there exists a knot ${\hat K}$ that is $n$-adjacent to the unknot and it
admits an \nt \ $L_n$ such that the interior of
$M_n:=M({\hat K}, L_n)$ admits a complete
hyperbolic structure with finite volume. Furthermore,
we have  ${\rm vol}(M_n)>n\,v_3 $, where
$v_3\,(\approx 1.01494)$ is the volume of a regular hyperbolic ideal tetrahedron.
\end{theo}

\proof To complete the proof we need the following claim:

\noindent {\bf Claim}\qua For every
$n>2$, there is a non-trivial  knot ${\hat K}$
that is not a satellite (i.e.\ its complement is atoroidal)
and it is $n$-adjacent to the unknot.

\noindent {\bf Proof of Claim}\qua 
By \cite{ak} there is a non-trivial knot $\bar K$
that is $n$-adjacent to the unknot.
Suppose $\bar K$ is a satellite knot.
By applying Theorem \ref{theo:satellite} inductively
we can find a pattern knot $\hat K$
of $\bar K$ such that i) $\hat K$ is not a non-trivial satellite;
and ii) 
$\hat K$ is $n$-adjacent to the unknot in a standard solid torus $V_1$.
By Lemma \ref{lem:nonisotopic},
$\hat K$ is not the unknot and the claim is proved.

To continue with the proof of the theorem, fix $n>2$ and let $\bar K$
be a knot as in the claim above.
Also let $L_n$ be an \nt \ for $\hat K$.
By Lemma \ref{lem:irreducible}, $M_n:=S^3\setminus \eta(L_n\cup {\bar K})$
is irreducible and $\partial$-incompressible.
We claim that $M_n$ is atoroidal. For, suppose that $M_n$ contains an essential torus.
By Lemma \ref{lem:zero}, $T$ must remain essential in the complement of $\hat K$. But this is impossible since $\hat K$ was chosen to be non-satellite. Hence $M_n$ is atoroidal and by Lemma \ref{lem:annuli} anannular.
Now Thurston's uniformization theorem  applies to conclude that $M_n$ is hyperbolic.
Since $\partial M_n$ has $n+1$ components, the interior of $M_n$
has $n+1$ cusps. By \cite{adams},
${\rm vol}(M_n)\geq ( n+1) v_3$ and thus
${\rm vol}(M_n)> n v_3$ as desired. \endproof

As we proved earlier,
Theorem \ref{theo:main} follows from Theorem \ref{theo:hyperbolic}.
Note that Corollary \ref{corol:Alexander} follows immediately from
Theorem \ref{theo:main}. We will finish the section with the proof
of Corollary \ref{corol:infinite}.

\begin{corol} For every $n, m\in {\bf N}$,
there exists a hyperbolic knot $K$ with trivial finite type invariants of orders $\leq m$
and such that
${\rm vol }(S^3\setminus K)>n$.
\end{corol}
\proof  Fix $ m\in {\bf N}$. By Theorem \ref{theo:main}, there is a hyperbolic
knot $K$ with ${\rm vol }(S^3\setminus K)>m$ and with trivial finite type invariants
of orders $\leq m$. Now this $K$ works for all $n$ with $n<m$.
To obtain a knot corresponding to $m,n$ for some $n>m$ just apply Theorem \ref{theo:main}
for this $n$. \endproof

\section {Concluding discussion}

\subsection{ Mahler measure of Alexander polynomials and volume} For a hyperbolic link $L$ with an unknotted component
$K$ let $L(q)$ denote the link obtained from $L\setminus K$
by a twist of order $q$ along $K$. Thurston's hyperbolic Dehn filling theorem implies
that as $q\longrightarrow \infty$, the volume of the complement of
$L(q)$  converges to the volume of the complement of $L$.
In \cite{sw}, D. Silver and S. Williams showed that the Mahler measure
of the Alexander polynomial of $L(q)$ exhibits an analogous behavior;
it converges to the Mahler measure
of the Alexander polynomial of $L$, as $q\longrightarrow \infty$.
This analogy prompted the question of whether the Mahler measure of
the Alexander polynomial of
a hyperbolic
knot is  related to the
volume of its complement \cite{sw}. In that spirit,
Corollary \ref{corol:Alexander} provides examples of non-alternating, hyperbolic knots with arbitrarily large
volume whose Alexander polynomial has trivial Mahler measure. In \cite{sstw},
Silver, Stoimenow and Williams construct examples of alternating  knots
with these properties. However they also show that, for alternating knots,
large
Mahler measure of the Alexander polynomial implies
large volume of the complement.

\subsection {The Jones polynomial and volume} Corollary \ref{corol:Alexander}
shows that the Alexander polynomial
of generic (non-alternating) hyperbolic knots is far from controlling the hyperbolic volume.
The determinant of a knot can also be
evaluated from its Jones polynomial. More specifically,
if 
$J_K(t):= \sum_{i=0}^r a_it^i$
denotes the Jones  polynomial of an alternating
knot $K$, then
${\rm det}(K):= |J_K(-1)|=\sum_{i=0}^r|a_i|$.
Thus, the inequalities in (4)
can also be expressed in terms of the Jones polynomial of $K$. 
In particular, the right hand side inequality becomes 
${\rm vol}(S^3\setminus K)\leq C(\sum_{i=0}^r|a_i|)$.
A recent result
O. Dasbach and X.-S. Lin \cite{dl}, which preceded and inspired the results in this
note,
gives a striking relation between the Jones polynomial
and the volume of alternating hyperbolic knots. More specifically,
the {\sl Volume-ish theorem} of \cite{dl} states that
the volume
of the complement
of a hyperbolic alternating knot
is bounded above and below by linear functions of absolute values of certain coefficients
of the Jones polynomial of the knot. Furthermore, the authors record
experimental data on the correlation between the coefficients of the Jones polynomial and the
volume of knots up to 14 crossings. For generic knots, this data suggests
a better
correlation than the one
between the coefficients of Alexander polynomial and volume.
In the view of the experimental evidence of \cite{dl},
and the examples constructed in this note, it becomes interesting to study the Jones polynomial of knots that are $n$-adjacent to the unknot,
for $n\gg0$.

\Addresses\recd
\end{document}